\documentclass[article,12pt]{amsart}

\usepackage{amsfonts}
\usepackage{amsmath}
\usepackage{mathrsfs}
\usepackage{graphicx}
\usepackage{color}
\usepackage{epsfig}
\usepackage{yfonts}
\usepackage{fancybox}
\usepackage{enumerate}
\usepackage{dsfont}

\numberwithin{equation}{section}
\theoremstyle{plain}
\newtheorem{thm}{Theorem}[section]
\newtheorem{rem}{Remark}[section]

\newcommand{\dE}{\mathbb{E}}

\newcommand{\dL}{\mathbb{L}}

\newcommand{\dZ}{\mathbb{Z}}
\newcommand{\dC}{\mathbb{C}}

\newcommand{\cF}{\mathcal{F}}

\newcommand{\cR}{\mathcal{R}}

\newcommand{\ind}{\mbox{1}\kern-.25em \mbox{I}}

\def\build#1_#2^#3{\mathrel{\mathop{\kern 0pt#1}\limits_{#2}^{#3}}}

\def\videbox{\mathbin{\vbox{\hrule\hbox{\vrule height1.4ex \kern.6em\vrule height1.4ex}\hrule}}}
\def\demend{\hfill $\videbox$\\}

\keywords{Elephant random walk, Martingales, strong law of large numbers, asymptotic normality}

\begin{document}
\title[Hypergeometric identities arising from the elephant random walk]
{Hypergeometric identities arising from the elephant random walk \vspace{1ex}}
\author{Bernard Bercu, Marie-Line Chabanol, Jean-Jacques Ruch}
\dedicatory{\normalsize University of Bordeaux, France}
\address{Universit\'e de Bordeaux, Institut de Math\'ematiques de Bordeaux,
UMR 5251, 351 Cours de la Lib\'eration, 33405 Talence cedex, France.}
\thanks{}

\begin{abstract}
A probabilistic approach is provided to establish new hypergeometric identities. It is based on the calculation
of moments of the limiting distribution of the position of the elephant random walk in the
superdiffusive regime.
\end{abstract}
\maketitle

\ \vspace{-7ex}
\section{Introduction}
\label{S-I}

Random walks with long-memory arose naturally in applied mathematics and statistical physics. One of them is the so-called 
elephant random walk (ERW) introduced by Sch\"utz and Trimper \cite{Schutz04} in the early 2000s, in order to investigate how long-range memory 
affects the random walk and induces a crossover from a diffusive to superdiffusive behavior. The ERW is a one-dimensional discrete-time random 
walk on $\dZ$, which has a complete memory of its entire history. It was referred to as the ERW in allusion to the famous saying that elephants 
can always remember where they have been before. 
\ \vspace{2ex} \par
The ERW starts at the origin at time zero,
$S_0=0$. At time $n=1$, the elephant moves to the right with probability $q$ and to the left with probability
$1-q$ where $q$ lies between $0$ and $1$, which means that $S_1=X_1$ where $X_1$ has a Rademacher $\cR(q)$ distribution.
Then,  at any time $n+1 \geq 1$, the elephant chooses uniformly at random an integer $k$ among the previous
times $1,\ldots,n$ and it moves exactly in the same direction as that of time $k$ with probability $p$ or in the opposite direction
with probability $1-p$. In other words,
\begin{equation}
\label{STEPS}
   X_{n+1} = \left \{ \begin{array}{ccc}
    +X_{k} &\text{ with probability } & p, \vspace{2ex}\\
    -X_{k} &\text{ with probability } & 1-p.
   \end{array} \right.
\end{equation}
Therefore, the position of the ERW at time $n+1$ is given by
\begin{equation}
\label{POSERW}
S_{n+1}=S_{n}+X_{n+1}.
\end{equation}
The asymptotic behavior of the ERW is closely related to the value of the probability $p \in [0,1]$ called the memory of the ERW. 
The elephant random walk is said to be subcritical if $0 \leq p < 3/4$, critical if $p=3/4$ and supercritical if $3/4< p \leq 1$. In this paper, we shall focus
our attention on the supercritical case where the memory parameter $3/4< p \leq 1$. It has been shown recently by three different approaches \cite{Baur16},
\cite{Bercu18}, \cite{Coletti17} that
\begin{equation}
\label{ASCVG}
 \lim_{n \rightarrow \infty} \frac{S_n}{n^{2p-1}}=L \hspace{1cm} \text{a.s.}
\end{equation}
where $L$ is a non-degenerate random variable. However, to the best of our knowledge, the distribution
of $L$ is far from being known. The only thing we know about the distribution of $L$ is that it is not Gaussian \cite{Bercu18}
Baur and Bertoin \cite{Baur16} established a very successful connection 
between the ERW and generalized P\'olya urns  \cite{Chauvin11}, \cite{Pemantle07}. Thanks to this connection and two functional limit theorems 
for multitype branching processes due to Janson \cite{Janson04}, it is possible to compute all the moments of $L$, see also \cite{Chauvin11}. 
The martingale approach of Bercu \cite{Bercu18} allows us to compute the moments of $L$ in two different ways. 
On the one hand, it is proven in \cite{Bercu18} that the almost sure convergence \eqref{ASCVG} 
also holds in $\dL^4$, which means that
\begin{equation}
\label{L4CVG}
 \lim_{n \rightarrow \infty} \dE\Bigl[ \Bigl| \frac{S_n}{n^{2p-1}} -L \Bigr|^4 \Bigr]=0.
\end{equation}
Consequently, in order to compute the first four moments of $L$, it is only necessary to calculate 
the first four moments of $S_n$ and to make use of  Lebesgue's dominated convergence theorem.
On the other hand, denote 
\begin{equation}
\label{DEFL}
 L_n = \frac{a_n S_n}{\Gamma(a+1)}
\end{equation}
where $a=2p-1$ and the sequence $(a_n)$ is given by $a_1=1$ and for all $n\geq 2$,
\begin{equation}
\label{DEFAN}
a_n =  \prod_{k=1}^{n-1}\Bigl(\frac{k+2p-1}{k}\Bigr)^{-1} = \frac{\Gamma(n)\Gamma(a+1)}{ \Gamma(n+a)} .
\end{equation}
We already saw in \eqref{L4CVG} that
\begin{equation}
\lim_{n \rightarrow \infty}L_n=L \hspace{1cm} \text{a.s.}
\end{equation}
Therefore, our strategy is to calculate the moments of $L$ in two different ways. 
We were happily surprised to find that this approach is the source of new hypergeometric identities. 
For example, we deduce from the double calculation of $\dE[L^2]$ that for all $a>1/2$,
\begin{equation*}
\frac{1}{(2a-1)\Gamma(2a)}=\frac{2a}{(2a-1)\Gamma(a+1)^2}\left(1 -\Bigl(\frac{a}{a+1}\Bigr)^2{}_{3}F_2 \Bigl( \begin{matrix}
{1,1,2a+1}\\
{a+2,a+2}\end{matrix} \Bigl|
1 \Bigr)\right)
\end{equation*}
leading to
\begin{equation*}
a^2 {}_{3}F_2 \Bigl( \begin{matrix}
{1,1,2a+1}\\
{a+2,a+2}\end{matrix} \Bigl|
1 \Bigr)=(a+1)^2-\frac{(\Gamma(a+2))^2}{\Gamma(2a+1)}
\end{equation*}
where $\!{}_{p}F_q$ stands for the hypergeometric function defined for all $z \in \dC$ by
\begin{equation}
\label{HGpq}
{}_{p}F_q \Bigl( \begin{matrix}
{a_1,\ldots,a_p}\\
{b_1,\ldots,b_q}\end{matrix} \Bigl|
{\displaystyle z}\Bigr)
=\sum_{n=0}^{\infty}
\frac{(a_1)_n\,\cdots\,(a_p)_n}
{(b_1)_n\,\cdots\,(b_q)_n\, n!} z^n.
\end{equation}
Here, $(a)_n$ is the rising Pochhammer symbol given by $(a)_0=1$ and, for all $n\geq 1$,
$(a)_n=a(a+1)\cdots(a+n-1)$. Via the double calculation of $\dE[L^d]$ for $d=2,3,4$,
we will obtain new hypergeometric identities associated to the function ${}_{d+1}F_d$.
\ \vspace{2ex} \par
A wide range of literature is available on hypergeometric identities.
We refer the reader to the very complete survey of Zudilin \cite{Zudilin} as well as to the more recent contributions
of Grondin \cite{Grondin} and Milgram \cite{Milgram}.
Our work is in the spirit of the classical Watson theorem, see the summation formulas in Lavoie \cite{Lavoie} and
Maier \cite{Maier}.
Krattenthaler and Srivinasa Rao \cite{Kratten} develop integral methods to prove many identities.
To our best knowledge, our identities are not one of them. 
\ \vspace{2ex} \par
Our main tool is the very successful method of telescoping sums that has been widely used by Zeilberger \cite{Zeilberger} 
and implemented in computer algebra systems in \cite{PWZ}. We would like to point out
that our identities were not directly recognized by Maple or Mathematica.
\ \vspace{2ex} \par
The paper is organized as follows. 
The main results of the paper are given in Section \ref{S-MR}. In Section \ref{S-M2}, we compute in two different ways the first order moments 
of the random variable $L$. It allows us to prove our new hypergeometric identities for the functions ${}_{3}F_2$, ${}_{4}F_3$ and ${}_{5}F_4$ 
evaluated at $z=1$.
Section \ref{S-Md} is devoted to the proof of the general case involving the higher order moments of $L$.
A short conclusion and perspectives are postponed to Section \ref{S-CP}.

\section{Main results}
\label{S-MR}

In order to simplify the presentation of our main results, denote $a=2p-1$ and $b=2q-1$.  We recall from \cite{Bercu18} that the first four moments of $L$ are given by 
\begin{eqnarray}
\label{MOML1}
\dE[L]  &=& \frac{b}{\Gamma(a+1)},\\
\label{MOML2}
\dE[L^2] &=& \frac{1}{(2a-1)\Gamma(2a)}, \\
\label{MOML3}
\dE[L^3] &=& \frac{b(a+1)}{a(2a-1)\Gamma(3a)}, \\
\label{MOML4}
\dE[L^4] &=& \frac{6(2a^2+2a-1)}{(4a-1)(2a-1)^2\Gamma(4a)}.
\end{eqnarray}

\noindent
By the calculation of the moments of $L$ in two different ways, we obtain the following  hypergeometric identities.
\begin{thm}
\label{T-NHI-MOM}
For any $a>1/2$, we have
\begin{equation}
\label{T-NHI1}
\frac{1}{(2a-1)\Gamma(2a)}= 
\frac{2a}{(2a-1)\Gamma(a+1)^2}\left(1 -\Bigl(\frac{a}{a+1}\Bigr)^2{}_{3}F_2 \biggl( \begin{matrix}
{1,1,2a+1} \\
{a+2,a+2}\end{matrix} \Bigl|
1 \Bigr)\right).
\end{equation}
\begin{eqnarray}
\frac{(a+1)}{a(2a-1)\Gamma(3a)} & = & \frac{3(a+1)}{(2a-1)\Gamma(a+1)^3}\left(1- \frac{3a^2}{(a+1)^3}
{}_{4}F_3 \Bigl( \begin{matrix}
{1,1,2,3a+1} \nonumber\\
{a+2,a+2,a+2}\end{matrix} \Bigl|
1 \Bigr)\right. \\
& - & 
\left.\frac{a^3}{(a+1)^3}
{}_{4}F_3 \Bigl( \begin{matrix}
{1,1,1,3a+1} \\
{a+2,a+2,a+2}\end{matrix} \Bigl|
1 \Bigr)\right).
\label{T-NHI2}
\end{eqnarray}
\begin{eqnarray}
\frac{6(2a^2+2a-1)}{(4a-1)(2a-1)^2\Gamma(4a)}& = & \frac{24a(2a^2+2a-1)}{(4a-1)(2a-1)^2\Gamma(a+1)^4}\biggl(1 
\nonumber \\
& - & 
\left.
\frac{6a^2}{(a+1)^4}
{}_{5}F_4 \Bigl( \begin{matrix}
{1,1,2,2,4a+1} \nonumber\\
{a+2,a+2,a+2,a+2}\end{matrix} \Bigl|
1 \Bigr)\right. \\
& - & 
\left.\frac{4a^3}{(a+1)^4}
{}_{5}F_4 \Bigl( \begin{matrix}
{1,1,1,2,4a+1} \nonumber\\
{a+2,a+2,a+2,a+2}\end{matrix} \Bigl|
1 \Bigr)\right. \\
& - &
\left.
\frac{a^4}{(a+1)^4}
{}_{5}F_4 \Bigl( \begin{matrix}
{1,1,1,1,4a+1} \\
{a+2,a+2,a+2,a+2}\end{matrix} \Bigl|
1 \Bigr)
\right).
\label{T-NHI3}
\end{eqnarray}
\end{thm}

\begin{rem}
Another identity coming from the calculation of $\dE[L^2]$ is given, for any $a>1/2$, by
\begin{eqnarray}
\label{T-NHI4}
\frac{1}{(2a-1)\Gamma(2a)} & = & \frac{1}{\Gamma(a+1)^2}
\!{\ }_{3}F_2 \Bigl( \begin{matrix}
{1,1,1} \\
{a+1,a+1}\end{matrix} \Bigl|
1 \Bigr)\\
   +\frac{ a^2 }{(2a-1)\Gamma(a+2)^2} &  & \hspace{-1cm} \left(
{ }_{3}F_2 \Bigl( \begin{matrix}
{1,1,2}\\
{a+2,a+2}\end{matrix} \Bigl|
1 \Bigr)
- 2a\, 
{ }_{3}F_2 \Bigl( \begin{matrix}
{1,1,2a+1}\\
{a+2,a+2}\end{matrix} \Bigl|
1 \Bigr)\right). \nonumber
\end{eqnarray}
\end{rem}

\noindent
A more general hypergeometric identity is as follows.

\begin{thm}
\label{T-NHI-GEN}
For any $a>1/2$ and for all $d \geq 2$, we have
\begin{equation}
\label{T-NHIGEN}
\sum_{k=0}^{d-2} \binom{k}{d}a^{d-k}
{}_{d+1}F_d \Bigl( \begin{matrix}
{1_{d-k},2_k, ad+1} \\
{a+2,a+2,\ldots,a+2}\end{matrix} \Bigl|
1 \Bigr)
=(a+1)^d - \frac{\Gamma(a+2)^d}{\Gamma(ad+1)}
\end{equation}
where $1_{d-k}=1,\dots,1$ of length $d-k$ and $2_{k}=2,\dots , 2$ of length $k$.
\end{thm}

\section{Proofs of the first order moment results}
\label{S-M2}

We shall now proceed to the proof of Theorem \ref{T-NHI-MOM}. Our proof relies on the Lebesgue's dominated convergence theorem
together with a telescoping argument. We already saw that the sequence
$(L_n)$ satisfies
$$
\sup_{n \geq 1} \dE[L_n^4] < \infty.
$$
Hence, the dominated convergence theorem implies that for any integer $d \leq 4$,
\begin{equation}
\label{LEBLN}
\lim_{n\rightarrow \infty} \dE[L_n^d]=\dE[L^d].
\end{equation}
where the limiting random variable $L$ is given by \eqref{DEFL}. Therefore, we have 
\begin{eqnarray}
 \dE[L^d]\Gamma(a+1)^d- \dE[X_1^d] & = & \lim_{n\rightarrow \infty} \dE[L_n^d]\Gamma(a+1)^d- \dE[X_1^d],\nonumber \\
 & & \nonumber \\
 & = & 
 \lim_{n\rightarrow \infty} \dE[(a_nS_n)^d]- \dE[(a_1S_1)^d]\nonumber\\
 & = & \lim_{n\rightarrow \infty} \sum_{k=2}^n(\dE[(a_kS_k)^d]- \dE[(a_{k-1}S_{k-1})^d])\nonumber \\
 & = & \sum_{k=2}^\infty\dE[(a_kS_k)^d]- \dE[(a_{k-1}S_{k-1})^d].
 \label{TEL}
\end{eqnarray}
We shall now proceed to the proof of \eqref{T-NHI1} by the computation of the second moment of the random variable $L$.
We have from identity (B.19) in \cite{Bercu18} that
\begin{equation}
\label{MOMSN2}
\dE[S_n^2]=\frac{n}{2a -1}\left( \frac{\Gamma(n+2a)}{\Gamma(n+1)\Gamma(2a)} -1\right).
\end{equation}
Consequently, we obtain from \eqref{DEFAN} together with \eqref{TEL} and \eqref{MOMSN2} that 
\begin{align*}
 \dE[L^2] & \Gamma(a+1)^2 -1  =\sum_{n=2}^\infty (a_n^2\dE[S_n^2]-a_{n-1}^2\dE[S_{n-1}^2])\\
 = &\frac{\Gamma(a+1)^2 }{2a-1} \sum_{n=2}^\infty\left[\frac{1}{\Gamma(2a)}\left(\frac{\Gamma(n)\Gamma(n+2a)}{\Gamma(n+a)^2}-\frac{\Gamma(n-1)\Gamma(n-1+2a)}{\Gamma(n-1+a)^2}\right) \right.\\
 &\hskip2.5cm -\left. \frac{\Gamma(n)\Gamma(n+1)}{\Gamma(n+a)^2}+\frac{\Gamma(n-1)\Gamma(n)}{\Gamma(n-1+a)^2}\right].
 \end{align*}
However, we have by the cancellation technique that
\begin{equation*}
\sum_{n=2}^\infty\left[ \frac{\Gamma(n)\Gamma(n+1)}{\Gamma(n+a)^2}-\frac{\Gamma(n-1)\Gamma(n)}{\Gamma(n-1+a)^2}\right]=-\frac{\Gamma(1)\Gamma(2)}{\Gamma(a+1)^2}=-\frac{1}{\Gamma(a+1)^2}
\end{equation*}
Therefore, we deduce from the previous calculation that
\begin{align*}
 \dE[L^2]-\frac{1}{\Gamma(a+1)^2} 
  = &\frac{-a^2}{(2a-1)\Gamma(2a)}\sum_{n=0}^\infty\frac{\Gamma(n+1)\Gamma(n+2a+1)}{\Gamma(n+2+a)^2} +\frac{1}{(2a-1)\Gamma(a+1)^2}\\
   = &\frac{-a^2}{(2a-1)\Gamma(2a)}
   {}_{3}F_2 \Bigl(\begin{matrix}
{1,1,2a+1} \\
{a+2,a+2}\end{matrix} \Bigl|
1 \Bigr)\frac{\Gamma(2a+1)}{\Gamma(a+2)^2}+\frac{1}{(2a-1)\Gamma(a+1)^2}
 \end{align*}
which clearly leads to \eqref{T-NHI1}. 
We are now in position to prove \eqref{T-NHI2} by the calculation of the third moment of $L$. It follows from identity (B.21) in \cite{Bercu18} that
 \begin{equation}
\dE[S_n^3] = \frac{b}{(2a -1)\Gamma(n)}\left(
\frac{3(a +1)\Gamma(n+3a)}{\Gamma(3a + 1)}
-\frac{\Gamma(n+a)}{\Gamma(a+1)}(3n+a +1)
\right).
\label{MOMSN3}
\end{equation}
Hence, we have from \eqref{TEL} and \eqref{MOMSN3} that 
\begin{align*}
 \dE[L^3] & \Gamma(a+1)^3 -b =\sum_{n=2}^\infty (a_n^3\dE[S_n^3]-a_{n-1}^3\dE[S_{n-1}^3])\\
= &\frac{b\Gamma(a\!+\!1)^3}{2a-1} \sum_{n=2}^\infty\left[\frac{3(a+1)}{\Gamma(3a+1)}\!\left(\frac{\Gamma(n)^3}{\Gamma(n\!+\!a)^3}\frac{\Gamma(n\!+\!3a)}{\Gamma(n)}\!-\!\frac{\Gamma(n\!-\!1)^3}{\Gamma(n\!-\!1\!+\!a)^3}\frac{\Gamma(n\!-\!1\!+\!3a)}{\Gamma(n-1)}\right)\right.\nonumber\\
 & -\frac{3}{\Gamma(a+1)}\left(\frac{\Gamma(n)^3}{\Gamma(n+a)^3}\frac{\Gamma(n+1+a)}{\Gamma(n)}-\frac{\Gamma(n-1)^3}{\Gamma(n-1+a)^3}\frac{\Gamma(n+a)}{\Gamma(n-1)}\right)\nonumber\\
 & \left.+\frac{2a-1}{\Gamma(a+1)}\left(\frac{\Gamma(n)^3}{\Gamma(n+a)^3}\frac{\Gamma(n+a)}{\Gamma(n)}-\frac{\Gamma(n-1)^3}{\Gamma(n-1+a)^3}\frac{\Gamma(n-1+a)}{\Gamma(n-1)}\right)\right].
 \end{align*}
On one hand, 
\begin{align*}
 & \sum_{n=2}^\infty\frac{3(a+1)}{\Gamma(3a+1)}\left(\frac{\Gamma(n)^3}{\Gamma(n+a)^3}\frac{\Gamma(n+3a)}{\Gamma(n)}-\frac{\Gamma(n-1)^3}{\Gamma(n-1+a)^3}\frac{\Gamma(n-1+3a)}{\Gamma(n-1)}\right)\\
& =\frac{3(a+1)}{\Gamma(3a+1)}\sum_{n=2}^\infty\frac{\Gamma(n-1)^2\Gamma(n-1+3a)}{\Gamma(n+a)^3}\left(-3(n-1)a^2-a^3\right)\\
 & =\frac{3(a+1)}{\Gamma(3a\!+\!1)}\!\left[-3a^2\sum_{n=0}^\infty\frac{\Gamma(n\!+\!1)\Gamma(n\!+\!2)\Gamma(n\!+\!1\!+\!3a)}{\Gamma(n+2+a)^3}\!-\!a^3\sum_{n=0}^\infty\frac{\Gamma(n\!+\!1)^2\Gamma(n\!+\!1\!+\!3a)}{\Gamma(n+2+a)^3}\right]\\
 & =\frac{3(a+1)}{\Gamma(a+2)^3}\left[-3a^2
{}_{4}F_3 \Bigl( \begin{matrix}
{1,1,2,3a+1} \nonumber\\
{a+2,a+2,a+2}\end{matrix} \Bigl|
1 \Bigr) 
 -  
a^3
{}_{4}F_3 \Bigl( \begin{matrix}
{1,1,1,3a+1} \\
{a+2,a+2,a+2}\end{matrix} \Bigl|
1 \Bigr)
\right].
  \end{align*}
On the other hand,  once again by the cancellation argument, we clearly have  
\begin{equation*}
\sum_{n=2}^\infty\left(\frac{\Gamma(n)^3}{\Gamma(n+a)^3}\frac{\Gamma(n+1+a)}{\Gamma(n)}-\frac{\Gamma(n-1)^3}{\Gamma(n-1+a)^3}\frac{\Gamma(n+a)}{\Gamma(n-1)}\right)=-\frac{\Gamma(a+2)}{\Gamma(a+1)^3}=-\frac{a+1}{\Gamma(a+1)^2}
\end{equation*}
and
\begin{equation*}
\sum_{n=2}^\infty \left(\frac{\Gamma(n)^3}{\Gamma(n+a)^3}\frac{\Gamma(n+a)}{\Gamma(n)}-\frac{\Gamma(n-1)^3}{\Gamma(n-1+a)^3}\frac{\Gamma(n-1+a)}{\Gamma(n-1)}\right)=-\frac{1}{\Gamma(a+1)^2}.
\end{equation*}
Putting together those three contributions, we obtain that
\begin{align}
 \dE[L^3] & -\frac{b}{\Gamma(a+1)^3} = \frac{3b(a+1) }{(2a-1)\Gamma(a+2)^3}\left(-3a^2{}_{4}F_3 \Bigl( \begin{matrix}
{1,1,2,3a+1} \nonumber\\
{a+2,a+2,a+2}\end{matrix} \Bigl|
1 \Bigr) \right. \\
& \left. -a^3
{}_{4}F_3 \Bigl( \begin{matrix}
{1,1,1,3a+1} \\
{a+2,a+2,a+2}\end{matrix} \Bigl|
1 \Bigr)
\right)+\frac{b}{2a-1}\left(\frac{3(a+1)}{\Gamma(a+1)^3}-\frac{(2a-1)}{\Gamma(a+1)^3}\right)\nonumber
 \end{align}  
which can be rewritten as
 \begin{align}
\hskip0.5cm \dE[L^3] & = \frac{3b(a+1) }{(2a-1)\Gamma(a+2)^3}\left(-3a^2{}_{4}F_3 \Bigl( \begin{matrix}
{1,1,2,3a+1} \\
{a+2,a+2,a+2}\end{matrix} \Bigl|
1 \Bigr)\right.\nonumber \\\
& \left. -a^3
{}_{4}F_3 \Bigl( \begin{matrix}
{1,1,1,3a+1} \\
{a+2,a+2,a+2}\end{matrix} \Bigl|
1 \Bigr)
\right)+\frac{3b(a+1)}{(2a-1)\Gamma(a+1)^3} \nonumber
\end{align}
which is exactly what we wanted to prove.
Now we are going to establish \eqref{T-NHI3} by the calculation of the fourth moment of $L$. It follows from identity (B.22)  in \cite{Bercu18} that
\begin{align}
\dE[S_n^4] = &\frac{1}{(2a-1)^2\Gamma(n)} \left(\frac{\Gamma(n+4a)}{(4a-1)\Gamma(4a)}6(2a^2+2a-1))-2(3n+2+2a)\frac{\Gamma(n+2a)}{\Gamma(2a)}\right.\nonumber\\
 &\left.+\frac{4a^2-12a+5}{(4a-1)}\Gamma(n+1)+3\Gamma(n+2)\right).
 \label{MOMSN4}
\end{align}
Using \eqref{TEL}, we get  that 
\begin{equation*}
\dE[L^4] \Gamma(a+1)^4-1= \sum_{n=2}^\infty (a_n^4\dE[S_n^4]-a_{n-1}^4\dE[S_{n-1}^4])
\end{equation*}
\begin{align}
 & = \frac{\Gamma(a\!+\!1)^4}{(2a\!-\!1)^2}\sum_{n=2}^\infty\left[\frac{6(2a^2\!+\!2a\!-\!1)}{(4a-1)\Gamma(4a)}\left(\!\frac{\Gamma(n)^4}{\Gamma(n\!+\!a)^4}\frac{\Gamma(n\!+\!4a)}{\Gamma(n)} \!-\!\frac{\Gamma(n\!-\!1)^4}{\Gamma(n\!-\!1\!+\!a)^4}\frac{\Gamma(n\!-\!1\!+\!4a)}{\Gamma(n-1)}\!\right)\right.\nonumber\\
 &\hskip0.5cm -\frac{6}{\Gamma(2a)}\left(\frac{\Gamma(n)^4}{\Gamma(n+a)^4}\frac{\Gamma(n+1+2a)}{\Gamma(n)}-\frac{\Gamma(n-1)^4}{\Gamma(n-1+a)^4}\frac{\Gamma(n+2a)}{\Gamma(n-1)}\right)\nonumber\\
 &\hskip0.5cm  \left.+\frac{4(2a-1)}{\Gamma(2a)}\left(\frac{\Gamma(n)^4}{\Gamma(n+a)^4}\frac{\Gamma(n+2a)}{\Gamma(n)}-\frac{\Gamma(n-1)^4}{\Gamma(n-1+a)^4}\frac{\Gamma(n-1+2a)}{\Gamma(n-1)}\right)\right.\nonumber\\
 &\hskip0.5cm  \left.+\frac{4a^2-12a+5}{(4a-1)}\left(\frac{\Gamma(n)^4}{\Gamma(n+a)^4}\frac{\Gamma(n+1)}{\Gamma(n)}-\frac{\Gamma(n-1)^4}{\Gamma(n-1+a)^4}\frac{\Gamma(n)}{\Gamma(n-1)}\right)\right.\nonumber\\
 &\hskip0.5cm  \left.+3\left(\frac{\Gamma(n)^4}{\Gamma(n+a)^4}\frac{\Gamma(n+2)}{\Gamma(n)}-\frac{\Gamma(n-1)^4}{\Gamma(n-1+a)^4}\frac{\Gamma(n+1)}{\Gamma(n-1)}\right)\right].\nonumber
 \end{align}
For the first part of the series, we have  
\begin{align*}
& \!\sum_{n=2}^\infty\!\frac{6(2a^2\!+\!2a\!-\!1)}{(4a-1)\Gamma(4a)}\left(\frac{\Gamma(n)^4}{\Gamma(n\!+\!a)^4}\frac{\Gamma(n\!+\!4a)}{\Gamma(n)}\!-\!\frac{\Gamma(n\!-\!1)^4}{\Gamma(n\!-\!1\!+\!a)^4}\frac{\Gamma(n\!-\!1\!+\!4a)}{\Gamma(n-1)}\!\right)\\
& =\frac{6(2a^2+2a-1)}{(4a-1)\Gamma(4a)}\left[-6a^2\sum_{n=0}^\infty\frac{\Gamma(n+1)\Gamma(n+2)^2\Gamma(n+1+4a)}{\Gamma(n+2+a)^4}\right.\\ 
&\left.-4a^3\sum_{n=0}^\infty\frac{\Gamma(n+1)^2\Gamma(n+2)\Gamma(n+1+4a)}{\Gamma(n+2+a)^4} -a^4 \sum_{n=0}^\infty\frac{\Gamma(n+1)^3\Gamma(n+1+4a)}{\Gamma(n+2+a)^4}\right]\\
 & =\frac{24(2a^2+2a-1)a}{(4a-1)\Gamma(a+2)^4}\left[-6a^2{}_{5}F_4 \Bigl( \begin{matrix}
{1,1,2,2,4a+1} \nonumber\\
{a+2,a+2,a+2,a+2}\end{matrix} \Bigl|
1 \Bigr) \right. \\
 &  \left.
 -4a^3 {}_{5}F_4 \Bigl( \begin{matrix}
{1,1,1,2,4a+1} \nonumber\\
{a+2,a+2,a+2,a+2}\end{matrix} \Bigl|
1 \Bigr)
 -a^4 {}_{5}F_4 \Bigl( \begin{matrix}
{1,1,1,1,4a+1} \\
{a+2,a+2,a+2,a+2}\end{matrix} \Bigl|
1 \Bigr) \right].
  \end{align*}
Using a well-known cancellation technique on the terms of the second, the third, the fourth and the fifth part of the series, 
we obtain that
\begin{align*}
& \sum_{n=2}^\infty\left(-\frac{\Gamma(n)^4}{\Gamma(n+a)^4}\frac{\Gamma(n+1+2a)}{\Gamma(n)}+\frac{\Gamma(n-1)^4}{\Gamma(n-1+a)^4}\frac{\Gamma(n+2a)}{\Gamma(n-1)}\right)
=\frac{\Gamma(2a+2)}{\Gamma(a+1)^4},\\
& \sum_{n=2}^\infty\left(\frac{\Gamma(n)^4}{\Gamma(n+a)^4}\frac{\Gamma(n+2a)}{\Gamma(n)}-\frac{\Gamma(n-1)^4}{\Gamma(n-1+a)^4}\frac{\Gamma(n-1+2a)}{\Gamma(n-1)}\right)=-\frac{\Gamma(2a+1)}{\Gamma(a+1)^4},\\
& \sum_{n=2}^\infty\left(\frac{\Gamma(n)^4}{\Gamma(n+a)^4}\frac{\Gamma(n+1)}{\Gamma(n)}-\frac{\Gamma(n-1)^4}{\Gamma(n-1+a)^4}\frac{\Gamma(n)}{\Gamma(n-1)}\right)
=-\frac{1}{\Gamma(a+1)^4},\\
& \sum_{n=2}^\infty\left(\frac{\Gamma(n)^4}{\Gamma(n+a)^4}\frac{\Gamma(n+2)}{\Gamma(n)}-\frac{\Gamma(n-1)^4}{\Gamma(n-1+a)^4}\frac{\Gamma(n+1)}{\Gamma(n-1)}\right)
=-\frac{2}{\Gamma(1+a)^4}.
\end{align*} 
Putting together all these contributions, we find that
\begin{eqnarray}
\dE[L^4] & = & \frac{24a(2a^2+2a-1)}{(4a-1)(2a-1)^2\Gamma(a+1)^4}\biggl(1 -\frac{6a^2}{(a+1)^4}
{}_{5}F_4 \Bigl( \begin{matrix}
{1,1,2,2,4a+1} \nonumber\\
{a+2,a+2,a+2,a+2}\end{matrix} \Bigl|
1 \Bigr)\\
& & -
\frac{4a^3}{(a+1)^4}
{}_{5}F_4 \Bigl( \begin{matrix}
{1,1,1,2,4a+1} \\
{a+2,a+2,a+2,a+2}\end{matrix} \Bigl|
1 \Bigr)\nonumber \\
&  & - \frac{a^4}{(a+1)^4}
{}_{5}F_4 \Bigl( \begin{matrix}
{1,1,1,1,4a+1} \\
{a+2,a+2,a+2,a+2}\end{matrix} \Bigl|
1 \Bigr) \biggl).\nonumber
\end{eqnarray}
In all our computations, $a=2p-1$ was a real number in $]1/2,1]$ but 
the hypergeometric identitities remain valid for any $a>1/2$ by analytic continuation.
\demend

\section{Proof of the higher order moment results}
\label{S-Md}
The computation of the first order moments of the random variable $L$ has led us to different identities involving the functions 
${}_{3}F_2$, ${}_{4}F_3$ and ${}_{5}F_4$ evaluated at $z=1$.
It is natural to expect some more general identity coming from the higher order moments of $L$.
This will be the core of the proof of Theorem \ref{T-NHI-GEN}.
\ \vspace{2ex}\par
First of all, in order to understand well how the elephant moves, denote by $(\cF_n)$ the increasing sequence of $\sigma$-algebras,
$\cF_n=\sigma(X_1,\ldots,X_n)$. It is straightforward to see from \eqref{STEPS} that the
conditional distribution $\mathcal{L}(X_{n+1}|\cF_n)=\mathcal{R}(p_n)$, which means that
\begin{equation}
\label{NSTEPS}
   X_{n+1} = \left \{ \begin{array}{ccc}
    1 &\text{ with probability } & p_n \vspace{2ex}\\
    -1 &\text{ with probability } & 1-p_n
   \end{array} \right.
\end{equation}
where
$$
p_n=\frac{1}{2}\Bigl(1+a\,\frac{S_n}{n}\Bigr).
$$
As a matter of fact, for any time $n \geq 1$, we have from \eqref{STEPS} that
\begin{equation*}
\dE[X_{n+1} | \cF_{n}] =a\frac{S_n}{n} \hspace{1cm}\mbox{a.s.}
\end{equation*}
Consequently, for all $d\ge 1$, we deduce $\dE[X_{n+1}^d | \cF_{n}]=1$ if $d$ is even, while if $d$ is odd
\begin{equation*}
\dE[X_{n+1}^d | \cF_{n}] = a\frac{S_n}{n} \hspace{1cm}\mbox{a.s.}
\end{equation*}
On the one hand, if $d$ is even, $d=2d'+2\ge 1$, we obtain from \eqref{POSERW} that 
\begin{eqnarray}
 \dE[S_{n+1}^{d} | \cF_{n}] & = &  \dE[(S_n+X_{n+1})^{d}| \cF_{n}]=\sum_{k=0}^{d}\binom{d}{k}S_n^{d-k} \dE[X_{n+1}^k | \cF_{n}], \nonumber \\
 & = & \sum_{k=0}^{d'+1}\binom{d}{2k}S_n^{d-2k} + \sum_{k=0}^{d'}\binom{d}{2k+1}S_n^{d-2k-1}a\frac{S_n}{n}, \nonumber \\
 & = &  1+ \sum_{k=0}^{d'} \left(\binom{d}{2k}+\frac{a}{n}\binom{d}{2k+1}\right)S_n^{d-2k}.
 \label{MCDSNE}
 \end{eqnarray}
On the other hand, if $d$ is odd, $d=2d'+1\ge 1$, we find from \eqref{POSERW} that 
\begin{eqnarray}
 \dE[S_{n+1}^{d} | \cF_{n}] & = &  \dE[(S_n+X_{n+1})^{d}| \cF_{n}]=\sum_{k=0}^{d}\binom{d}{k}S_n^{d-k} \dE[X_{n+1}^k | \cF_{n}], \nonumber\\
 & = & \sum_{k=0}^{d'}\binom{d}{2k}S_n^{d-2k} + \sum_{k=0}^{d'}\binom{d}{2k+1}S_n^{d-2k-1}a\frac{S_n}{n},\nonumber\\
 & = &  \sum_{k=0}^{d'} \left(\binom{d}{2k}+\frac{a}{n}\binom{d}{2k+1}\right)S_n^{d-2k}.
 \label{MCDSN0}
 \end{eqnarray}
Identities \eqref{MCDSNE} and \eqref{MCDSN0} can be rewritten in a single form as follows. For $d$ even, $d=2d'+2$
or for $d$ odd, $d=2d'+1$, we have
\begin{equation}
\label{CONDMOMSND} 
 \dE[S_{n+1}^{d} | \cF_{n}] = \sum_{k=0}^{d'} A_{k,d,n} S_n^{d-2k}+ 
 \mathds{1}_{d \text{ is even}}
 \end{equation} 
where
$$
A_{k,d,n}=\binom{d}{2k}+\frac{a}{n}\binom{d}{2k+1}
$$
and
$$
\mathds{1}_{d \text{ is even}}= 
\left \{ \begin{array}{cc}
    1 &\text{ if } d  \text{ is even, } \vspace{2ex}\\
    0 &\text{ if } d \text{ is odd. }
   \end{array} \right.
$$
By taking the expectation on both sides of \eqref{CONDMOMSND}, we deduce that
\begin{equation}
\label{MOMSND} 
 \dE[S_{n+1}^{d}] = \sum_{k=0}^{d'} A_{k,d,n} \dE[S_n^{d-2k}]+ 
 \mathds{1}_{d \text{ is even}}
 \end{equation} 
which leads for all $n \geq 1$, to
\begin{equation}
\label{INDUC}
Z_{n+1}=A_n Z_n + B_n
\end{equation}
where  $Z_n=\dE[S_n^{d}]$, $A_n=A_{0,d,n}$ and
$$
B_n= \sum_{i=1}^{d'} A_{i,d,n} \dE[S_n^{d-2i}]+ 
\mathds{1}_{d \text{ is even}}.
$$
It is not hard to see that \eqref{INDUC} implies that
\begin{equation}
\label{NINDUC}
Z_{n+1} = Z_1\prod_{k=1}^{n} A_k + \sum_{k=1}^{n} B_k \prod_{\ell=k+1}^{n} A_\ell.
\end{equation}
Consequently, we deduce from \eqref{MOMSND} and \eqref{NINDUC} that
\begin{equation*}
\dE[S_{n+1}^{d}]= \dE[X_{1}^{d}]\prod_{k=1}^nA_{0,d,k}+\mathds{1}_{d \text{ is even}}\sum_{k=1}^{n} \!\prod_{\ell=k+1}^n \! A_{0,d,\ell}+ 
\sum_{i=1}^{d'} \sum_{k=1}^{n} A_{i,d,k}\!\!\prod_{\ell=k+1}^n\!\!A_{0,d,\ell}\dE[S_k^{d-2i}].
\end{equation*}
Replacing $A_{k,d,n}$ in the above identity, we obtain that
\begin{align}
\label{CALCSND}
 \dE[S_{n}^{d}] & =\frac{\Gamma(n+ad)}{\Gamma(n)} \left(\frac{1}{\Gamma(ad+1)}\dE[X_{1}^{d}]+\mathds{1}_{d \text{ is even}}\sum_{k=1}^{n-1} \frac{\Gamma(k+1)}{\Gamma(k+1+ad)}\right. \\
  + \sum_{i=1}^{d'} \sum_{k=1}^{n-1} &\left. \frac{\Gamma(d+1)}{\Gamma(d\!-\!2i\!+\!1)\Gamma(2i\!+\!2)}\left[2i\left(1\!-\!\frac{a}{k}\right)\!+\!1\!+\!\frac{ad}{k}\right]
  \frac{\Gamma(k+1)}{\Gamma(k\!+\!ad\!+\!1)}\dE[S_k^{d-2i}]\right). \nonumber
 \end{align}
Furthermore, we have from \eqref{DEFAN} that
\begin{equation}
\lim_{n\rightarrow \infty} \frac{a_n^d \Gamma(n+ad)}{\Gamma(n)}=\Gamma(a+1)^d.
\label{LIMAND}
\end{equation}
It clearly follows from \eqref{CALCSND}, \eqref{LIMAND} and by induction on $d$ that for all $d \geq 2$ even,
\begin{equation}
\label{SUPSNDEVEN}
\sup_{n \geq 1} a_n^d \dE[S_n^d] < \infty.
\end{equation}
Hence, we deduce from \eqref{SUPSNDEVEN} together with H\"{o}lder's inequality for all $d \geq 1$ odd,
\begin{equation*}
\label{SUPSNDODD}
\sup_{n \geq 1} a_n^d \dE[|S_n|^d] < \infty,
\end{equation*}
which ensures that for all integer $d \geq 1$, the sequence $(L_n)$ is bounded in $\dL^d$, 
\begin{equation}
\label{SUPLND}
\sup_{n \geq 1} \dE[|L_n|^d] < \infty.
\end{equation}
Consequently, the Lebesgue's dominated convergence theorem implies that for all integer $d \geq 1$,
\begin{equation}
\label{LIMLD}
\lim_{n\rightarrow \infty} \dE[L_n^d]=\dE[L^d].
\end{equation}
We are now in position to prove \eqref{T-NHIGEN}. We already saw from \eqref{TEL} that
\begin{equation}
 \dE[L^d]\Gamma(a+1)^d- \dE[X_1^d]  =  \sum_{k=2}^\infty\dE[(a_kS_k)^d]- \dE[(a_{k-1}S_{k-1})^d].
\end{equation}
On the one hand, we obtain from \eqref{DEFAN} and \eqref{CALCSND} that
\begin{align*}
 & \frac{a_n^d \Gamma(n+ad)}{\Gamma(n)} -  \frac{a_{n-1}^d\Gamma(n-1+ad)}{\Gamma(n-1)},\\
 & =  \frac{\Gamma(n)^{d-1}\Gamma(a+1)^{d} }{\Gamma(n+a)^{d}} \Gamma(n+ad) -
 \frac{\Gamma(n-1)^{d-1}\Gamma(a+1)^{d}}{\Gamma(n-1+a)^{d}}\Gamma(n-1+ad), \\
& =  \frac{\Gamma(a+1)^d \Gamma(n-1)^{d-1}}{\Gamma(n+a)^d}\Gamma(n-1+ad)\Bigl((n-1)^{d-1}(n-1+ad)-(n-1+a)^d\Bigr),\\
& = - \frac{\Gamma(a+1)^d\Gamma(n-1)^{d-1}}{\Gamma(n+a)^d}\Gamma(n-1+ad)\sum_{k=0}^{d-2}\binom{d}{k}(n-1)^ka^{d-k}.
\end{align*} 
Now, by summing on $n$, we obtain that
\begin{eqnarray}
 & & \sum_{n=2}^{\infty}\Bigl(\frac{a_n^d\Gamma(n+ad)}{\Gamma(n)\Gamma(ad+1)} -  \frac{a_{n-1}^d\Gamma(n-1+ad)}{\Gamma(n-1)\Gamma(ad+1)}\Bigr)
\nonumber\\
 & = & -\frac{\Gamma(a+1)^d}{\Gamma(ad+1)}\sum_{k=0}^{d-2}\binom{d}{k}a^{d-k}\sum_{n=2}^{\infty}\frac{\Gamma(n-1)^{d-1-k}\Gamma(n)^{k}}{\Gamma(n+a)^d}\Gamma(n-1+ad), \nonumber\\
 & = & -\frac{1}{(a+1)^d}\sum_{k=0}^{d-2} \binom{k}{d}a^{d-k}
{}_{d+1}F_d \Bigl( \begin{matrix}
{1_{d-k},2_k, ad+1} \\
{a+2,a+2,\ldots,a+2}\end{matrix} \Bigl|
1 \Bigr)
\label{GENHYP1}
\end{eqnarray} 
where $1_{d-k}=1,\dots,1$ of length $d-k$ and $2_{k}=2,\dots , 2$ of length $k$.
On the other hand, as we have a telescoping series,
\begin{equation}
\sum_{n=2}^{\infty} \Bigl( \frac{a_n^d \Gamma(n+ad)}{\Gamma(n)\Gamma(ad+1)} -  \frac{a_{n-1}^d\Gamma(n-1+ad)}{\Gamma(n-1)\Gamma(ad+1)} \Bigr) =\frac{\Gamma(a+1)^{d}}{ \Gamma(ad+1)}-1.
\label{GENHYP2}
\end{equation}
Finally, \eqref{GENHYP1} together with \eqref{GENHYP2} immediately lead to \eqref{T-NHIGEN}, which completes the proof
of Theorem \ref{T-NHI-GEN}.
\demend

\vspace{-3ex}
\section{Conclusion and perspectives}
\label{S-CP}

In this paper, we have computed the moments of the position of the ERW is two different ways.
It led us to find new hypergeometric identities. However, even if we know all the moments of $L$,
its distribution is far from being known. Baur and Bertoin \cite{Baur16} established a very successful connection 
between the ERW and generalized P\'olya urns. Unfortunately, the results of Chauvin and Pouyanne \cite{Chauvin11}
as well as the induction formulas of Janson \cite{Janson04} seem useless to go deeper into the knowledge
of the distribution of $L$. It would be very interesting to specify what is the exact distribution of $L$.
\ \vspace{2ex} \par
Another source of hypergeometric identities, coming from statistical physics, is the elephant random walk with stops \cite{Cressoni13}. 
The ERW by stops starts at the origin at time zero, $S_0=0$. At time $n=1$, the elephant moves to the right with probability $s$ and to the left with probability $1-s$ where $s$ lies between $0$ and $1$. Then,  at any time $n+1 \geq 1$, the elephant chooses uniformly at 
random an integer $k$ among the previous times $1,\ldots,n$ and we have
\begin{equation*}
   X_{n+1} = \left \{ \begin{array}{ccc}
    +X_{k} &\text{ with probability } & p, \vspace{1ex}\\
    -X_{k} &\text{ with probability } & q, \vspace{1ex}\\
    0 &\text{ with probability } & r,
   \end{array} \right.
\end{equation*}
where $p+q+r=1$. In the special case $r=0$, we find again the standard ERW.
By the same procedure as before, we obtain the following identity, which can be seen as a generalization of \eqref{T-NHI4}.
For any $a>0$ and $b>0$ such that $2a>b$,
\begin{eqnarray}
\label{T-NHI5}
\frac{1}{(2a-b)\Gamma(2a)} & = & \frac{1}{\Gamma(a+1)^2}
\!{\ }_{3}F_2 \Bigl( \begin{matrix}
{1,1,b} \\
{a+1,a+1}\end{matrix} \Bigl|
1 \Bigr) \\
   +\frac{ a^2 }{(2a-b)\Gamma(a+2)^2} &  & \hspace{-1cm} \left( b\,
{ }_{3}F_2 \Bigl( \begin{matrix}
{1,1,b+1}\\
{a+2,a+2}\end{matrix} \Bigl|
1 \Bigr)
- 2a\, 
{ }_{3}F_2 \Bigl( \begin{matrix}
{1,1,2a+1}\\
{a+2,a+2}\end{matrix} \Bigl|
1 \Bigr)\right). \nonumber
\end{eqnarray}

\section*{Acknowledgements}
The authors would like to thank Alin Bostan, Henri Cohen and Wadim Zudilin for fruitful discussions
on a preliminary version of this paper.
\bibliographystyle{acm}

\begin{thebibliography}{1}

\bibitem{Baur16}
{\sc Baur, E. and Bertoin, J.} 
\newblock Elephant Random Walks and their connection to P\'olya-type urns.
\newblock {\em Phys. Rev. E 94}, 052134 (2016).


\bibitem{Bercu18}
{\sc Bercu, B.}
\newblock A martingale approach for the elephant random walk.
\newblock {\em J. Phys. A: Math. Theor. 51}, 015201 (2018).


\bibitem{Chauvin11}
{\sc Chauvin, B., Pouyanne, N., Sahnoun, R.} 
\newblock Limit distributions for large P\'olya urns. 
\newblock {\em Ann. Appl. Probab. 21}, (2011), pp 1-32.

\bibitem{Coletti17}
{\sc Coletti, C. F., Gava, R., Sch\"utz, G. M.} 
\newblock Central limit theorem and related results for the elephant random walk.
\newblock {\em J. Math. Phys. 58}, 053303 (2017).

\bibitem{Cressoni13}
{\sc Cressoni, J. C., Viswanathan, G. M., Da Silva, M. A. A.} 
\newblock Exact solution of an anisotropic 2D random walk model with strong memory correlations.
\newblock {\em J. Phys. A: Math. Theor. 46}, 505002 (2013).


\bibitem{Da13}
{\sc Da Silva, M.~A.~A., Cressoni, J.~C., Sch\"utz, G.~M., Viswanathan, G.~M., Trimper, S.}
\newblock Non-Gaussian propagator for elephant random walks. 
\newblock {\em Phys. Rev. E 88}, 022115 (2013).


\bibitem{Grondin}
{\sc Grondin, F.},
\newblock On partial sums of generalized hypergeometric series.
Preprint, (2017).

\bibitem{Janson04}
{\sc Janson, }
\newblock Functional limit theorems for multitype branching processes and generalized P\'{o}lya urns. 
\newblock {\em Stochastic Process. Appl. 110}, (2004), pp. 177-245.


\bibitem{Kratten}
{\sc Krattenthaler, C. and Srinivasa Rao, K. }
\newblock Automatic generation of hypergeometric identities by the beta integral method
\newblock {\em Journal of Computational and Applied Mathematics 160}, (2003), pp. 159-173.

\bibitem{Kursten16}
{\sc K\"ursten, R.}
\newblock Random recursive trees and the elephant random walk. 
\newblock {\em Phys. Rev. E 93}, 032111 (2016).

\bibitem{Lavoie}
{\sc Lavoie, J. L.}
\newblock Some Summation Formulas for the Series $\!\!{\ }_3F_2(1)$.
\newblock {\em Mathematics of computation 49}, (1987), pp. 269-274.

\bibitem{Maier}
{\sc Maier, R.}
\newblock A generalization of Euler's hypergeometric transformation.
\newblock {\em Transaction of the AMS 358}, (2006), pp. 39-57. 

\bibitem{Milgram}
{\sc Milgram, M.}  
\newblock Variations on a hypergeometric theme.
\newblock {\em Journal of Classical Analysis 13}, (2018), pp. 1-43.

\bibitem{Pemantle07}
{\sc Pemantle, R.}
\newblock A survey of random processes with reinforcement. 
\newblock {\em Probability Surveys 4}, (2007), pp 1-79.

\bibitem{PWZ}
{\sc Petkovsek, M.,  Wilf, H. S.,  and Zeilberger, D.}
\newblock A=B
\newblock{A.K. Peters, CRC Press}, 1996.


\bibitem{Schutz04}
{\sc Sch\"utz, G.~M., and Trimper, S.}
\newblock Elephants can always remember: Exact long-range memory
effects in a non-Markovian random walk. 
\newblock {\em Phys. Rev. E 70}, 045101 (2004).

\bibitem{Zeilberger}
{\sc Zeilberger, D.} 
\newblock The method of creative telescoping.
\newblock {\em J. Symbolic Computation 11},
(1991), pp. 195-204. 

\bibitem{Zudilin}
{\sc Zudilin, W.}
\newblock Arithmetic hypergeometric series.
\newblock {\em Russian Math surveys, 66}, (2011), pp. 369-420.

\end{thebibliography}

\end{document}